\documentclass[11pt,a4paper,oneside,reqno]{amsart}
\usepackage[T1]{fontenc}
\usepackage{lmodern,amssymb}
\usepackage[a4paper,textwidth=162mm,top=23mm,bottom=23mm]{geometry}
\usepackage[hidelinks]{hyperref}
\hypersetup{
  pdftitle={The totient bound and extremal gaps in cyclotomic polynomials},
  pdfauthor={Qi-Yang Zheng}
}
\newtheorem{theorem}{Theorem}[section]
\newtheorem{conjecture}{Conjecture}[section]
\newtheorem{lemma}[theorem]{Lemma}
\newtheorem{corollary}[theorem]{Corollary}
\newtheorem{proposition}[theorem]{Proposition}
\theoremstyle{definition}
\newtheorem{example}[theorem]{Example}
\newcommand{\Ph}{\Phi_m}
\newcommand{\Q}{\mathbb Q}
\newcommand{\Z}{\mathbb Z}
\newcommand{\supp}{\operatorname{supp}}
\numberwithin{equation}{section}
\title[Extremal gaps in cyclotomic polynomials]{The totient bound and extremal gaps\\ in cyclotomic polynomials}
\author{\NoCaseChange{Qi-Yang Zheng}}
\date{}
\begin{document}
\begin{abstract}
Let \(m>1\) be odd and squarefree, and let \(p\nmid m\) be prime.
If \(p<m\), every gap in \(\Phi_{mp}\) has length strictly less than
\(\varphi(m)\). If \(p>m\), the maximum gap has length \(\varphi(m)\)
and occurs exactly \(2\lfloor p/m\rfloor\) times. We determine its
positions and endpoint coefficients. These results prove the equality
criterion of Ambrosino, Hong and Lee and the multiplicity formula of
Zhang. We also treat repeated prime factors, even cofactors and
extremal zero runs in periodic coefficient blocks.
\end{abstract}
\vspace*{0pt}
\maketitle
\vspace{-0.8\baselineskip}

\section{Introduction}
A primitive \(n\)th root of unity is a complex number of multiplicative
order \(n\). The \(n\)th cyclotomic polynomial is
\(\Phi_n(x)=\prod_{\zeta}(x-\zeta)\), where the product runs over all such roots. Its degree is Euler's totient
\(\varphi(n)=\#\{1\le k\le n:\gcd(k,n)=1\}\).
For \(F(x)=\sum_{k=0}^{D}c_kx^k\), write \(\supp F=\{k:c_k\ne0\}\).
A \emph{gap} is a pair \((a,b)\) of consecutive elements of \(\supp F\).
Its length is \(b-a\), so it contains \(b-a-1\) intermediate zero
coefficients. Write \(g(F)\) for the greatest gap length and \(N(F)\)
for the number of gaps of that length.

For odd primes \(q<p\), Hong, Lee, Lee and Park \cite{HLLP} proved
that \(g(\Phi_{qp})=q-1\). For example, with \(q=3\) and \(p=5\),
\[
   \Phi_{15}(x)=1-x+x^3-x^4+x^5-x^7+x^8.
\]
The consecutive exponent differences are \(1,2,1,1,2,1\), so
\(g(\Phi_{15})=2\). Camburu, Ciolan, Luca, Moree and Shparlinski
\cite{CCLMS} and Zhang \cite{ZB} independently showed that the
length \(q-1\) occurs \(2\lfloor p/q\rfloor\) times. In \(\Phi_{15}\)
the two occurrences are \((1,3)\) and \((5,7)\), in agreement with
\(2\lfloor5/3\rfloor=2\).

Two conjectures extend these results to products of distinct odd primes.
In both, write such a product as \(n=mp\), where \(p\) is its largest prime factor.

\begin{conjecture}[Ambrosino, Hong and Lee {\cite[Conjecture~7]{AHL}}]\label{conj:criterion}
The maximum gap of \(\Phi_n\) has length \(\varphi(m)\) if and only if \(p>m\).
\end{conjecture}

\begin{conjecture}[Zhang {\cite[Conjecture~0.1]{ZT}}]\label{conj:count}
If \(m>1\) and \(p>m\), then
\[
  g(\Phi_{mp})=\varphi(m),\qquad
  N(\Phi_{mp})=2\lfloor p/m\rfloor .
\]
\end{conjecture}

Al-Kateeb, Ambrosino, Hong and Lee \cite[Theorem~1]{AAHL} proved the length
formula for \(p>m\). Theorem~\ref{thm:main}(i) completes
Conjecture~\ref{conj:criterion} by proving strict inequality for \(p<m\).
Part~(ii) proves the count in Conjecture~\ref{conj:count} for any number
of prime factors and gives the positions and endpoint coefficients of
all maximum gaps.

\begin{theorem}\label{thm:main}
Let \(m>1\) be odd and squarefree, let \(p\) be a prime with \(p\nmid m\),
and put \(d=\varphi(m)\).

\smallskip\noindent\textbf{\textup{(i)} \(p<m\).}
The maximum gap length satisfies \(g(\Phi_{mp})<d\).

\smallskip\noindent\textbf{\textup{(ii)} \(p>m\).}
The maximum gap length is \(d\). Put \(D=d(p-1)\) and
\(Q=\lfloor p/m\rfloor\). There are exactly \(2Q\) gaps of length \(d\),
whose pairs of endpoints are
\begin{equation}\label{eq:pairs}
   (jm-d,jm) \qquad\text{and}\qquad (D-jm,D-jm+d),
   \qquad 1\le j\le Q .
\end{equation}
In the first pair, the coefficient of \(x^{jm-d}\) is \(-1\) and that
of \(x^{jm}\) is \(+1\). In the second, the coefficient of \(x^{D-jm}\)
is \(+1\) and that of \(x^{D-jm+d}\) is \(-1\). Every coefficient
strictly between the two endpoints of either pair is zero.
\end{theorem}

\noindent\textit{Example.}
Take \(m=15\), so \(d=8\). For \(p=13<15\), the maximum in \(\Phi_{195}\)
is \(6<8\) and occurs four times \cite[Theorem~0.1]{ZT}. For \(p=17>15\),
we have \(D=128\) and \(Q=1\). The two gaps of length \(8\) in \(\Phi_{255}\)
run from \(-x^7\) to \(+x^{15}\) and from \(+x^{113}\) to \(-x^{121}\).

When \(n=p\), Conjecture~\ref{conj:criterion} follows directly from
\(\Phi_p(x)=1+x+\cdots+x^{p-1}\) and \(\varphi(1)=1\).
For \(m>1\), Theorem~\ref{thm:main} allows any prime \(p\nmid m\).
Sanna \cite[\S6]{Sanna} records the proposed upper bound for \(p<m\).
Zhang \cite[Theorem~0.1]{ZT} determined the maximum and its multiplicity
for \(m=15\). Benoist and Perucca \cite{BP} proved the predicted count for
\(m=3q\), with \(q>3\) prime, \(p>3q\) and \(p\equiv\pm1,\pm2\pmod{3q}\).

The proof of Theorem~\ref{thm:main} uses the product
\(\Phi_m(x)\Phi_{mp}(x)=\Phi_m(x^p)\), which also occurs in Kaplan's
binary argument in \cite[\S2.1]{CCLMS}. Truncating it at the right
endpoint of a gap reduces the equality case to a question about a
partial sum of \(\Phi_m\).

\begin{corollary}[Odd cofactors]\label{cor:odd-radical}
Let \(m>1\) be odd and let \(p\nmid m\) be prime. Put
\(s=\operatorname{rad}(m)=\prod_{q\mid m,\ q\ \mathrm{prime}}q\),
\(d=\varphi(m)\) and \(D=d(p-1)\). Then
\[
  g(\Phi_{mp})\le d,\qquad
  g(\Phi_{mp})=d\quad\Longleftrightarrow\quad p>s .
\]
The gaps of length \(d\) are exactly
\[
   (jm-d,jm)\quad\text{and}\quad (D-jm,D-jm+d),
   \qquad 1\le j\le\lfloor p/s\rfloor .
\]
Their endpoint coefficients are \((-1,1)\) and \((1,-1)\), respectively.
Thus there are \(2\lfloor p/s\rfloor\) gaps of length \(d\).
\end{corollary}
\begin{proof}
For \(n>1\), put \(R=\operatorname{rad}(n)\) and \(v=n/R\).
In the M\"obius product for \(\Phi_n\), a divisor \(t\mid n\)
contributes precisely when \(n/t\mid R\), that is, when \(t=ve\)
with \(e\mid R\). This gives the classical radical reduction
\cite[\S2]{HLLP}
\begin{equation}\label{eq:radical}
   \Phi_n(x)=\prod_{e\mid R}(1-x^{ve})^{\mu(R/e)}=\Phi_R(x^v).
\end{equation}
Set \(h=m/s\). Since \(p\nmid m\), we have \(\operatorname{rad}(mp)=sp\)
and hence \(\Phi_{mp}(x)=\Phi_{sp}(x^h)\). Every gap \((a,b)\) of
\(\Phi_{sp}\) corresponds to the gap \((ha,hb)\) with the same endpoint
coefficients. Theorem~\ref{thm:main} applied to \(s\) now gives the
claim, since \(hs=m\) and \(h\varphi(s)=d\).
\end{proof}

For example, take \(m=9\) and \(p=5\). In \(\Phi_{45}(x)=\Phi_{15}(x^3)\)
the gaps \((1,3)\) and \((5,7)\) of \(\Phi_{15}\) become \((3,9)\)
and \((15,21)\). Both have length \(6=\varphi(9)\), as
\(\operatorname{rad}(9)=3<5<9\).

Section~\ref{sec:even} treats every even cofactor by the radical and sign
identities of \cite[\S2]{HLLP}. For squarefree even \(m>2\) and
\(p\nmid m\), the threshold is \(p>m/2\). The final section describes the
cyclic runs of zeros for both parities.

\section{Truncation and complementary partial sums}
Fix \(m>1\). Write \(\Ph(T)=\sum_{j=0}^{d}u_jT^j\), where
\(d=\varphi(m)\), and put \(P_i(T)=\sum_{j=0}^{i}u_jT^j\)
for \(0\le i<d\). We have \(u_0=u_d=1\).

\begin{proposition}\label{prop:trunc}
For every prime \(p\nmid m\) one has \(g(\Phi_{mp})\le d\).
Suppose that \((a,b)\) is a gap of length \(d\), and let
\(c=[x^a]\Phi_{mp}(x)\). Then for every primitive \(m\)th root \(\zeta\),
\begin{equation}\label{eq:unit}
   P_i(\zeta^p)=-c\zeta^b,\qquad
   i=\left\lfloor\frac{b-1}{p}\right\rfloor .
\end{equation}
\end{proposition}
\begin{proof}
Set \(F=\Phi_{mp}\). A root of unity of order \(n\) has \(p\)th power
of order \(n/\gcd(n,p)\). Since \(p\nmid m\), this order equals \(m\)
precisely when \(n=m\) or \(n=mp\). The roots of \(\Ph(x^p)\) are
simple because this polynomial divides \(x^{mp}-1\). Both factors
on the left below are monic, so comparison of roots and then degrees gives
\begin{equation}\label{eq:product}
   \Ph(x)F(x)=\Ph(x^p),\qquad D=\deg F=pd-d=d(p-1).
\end{equation}
For any gap \((a,b)\), let \(A(x)\) be the sum of the terms of \(F(x)\)
of degree less than \(b\). Since \(a\) is the last support exponent
before \(b\), we have \(\deg A=a\) and \(F\equiv A\pmod{x^b}\).
On the right of \eqref{eq:product}, the term \(u_jx^{pj}\) survives
truncation modulo \(x^b\) only in degrees \(pj<b\). The indices in this
range are exactly \(0\le j\le i=\lfloor(b-1)/p\rfloor\). Moreover,
\(b\le D=d(p-1)\), since \(b\) is an exponent of \(F\). Thus
\(ip\le b-1<dp\), which gives \(i<d\). Consequently
\begin{equation}\label{eq:trunc}
   \Ph(x)A(x)\equiv\sum_{j=0}^{i}u_jx^{pj}=P_i(x^p)\pmod{x^b}.
\end{equation}
If \(b-a>d\), the left side has degree \(d+a<b\), and the right side
has degree at most \(pi\le b-1\). Their difference is divisible by
\(x^b\), hence is zero. Evaluation at a primitive \(m\)th root
\(\zeta\) gives \(P_i(\zeta^p)=0\). But \(\zeta^p\) is again primitive
of order \(m\), so its minimal polynomial over \(\Q\) has degree \(d\).
This contradicts \(P_i(0)=1\) and \(\deg P_i\le i<d\). The gap bound
follows.

Now suppose \(b-a=d\). By \eqref{eq:trunc}, there is a polynomial
\(f(x)\) such that \(\Ph(x)A(x)-P_i(x^p)=x^b f(x)\). The difference
has degree exactly \(a+d=b\), since \(\deg P_i(x^p)<b\). Its leading
coefficient is \(c=[x^a]F(x)\), since \(\Ph\) is monic. The quotient \(f\) therefore has degree zero. Comparing the
coefficients of \(x^b\) gives \(f=c\), so
\[
   \Ph(x)A(x)=P_i(x^p)+cx^b .
\]
At \(x=\zeta\) the left side vanishes, giving \eqref{eq:unit}.
\end{proof}

\begin{lemma}\label{lem:prefix}
Suppose \(m\) is odd and squarefree. For a primitive \(m\)th root
\(\eta\) and \(1\le i\le d-2\), one has \(\lvert P_i(\eta)\rvert^2\notin\Q\).
\end{lemma}
\begin{proof}
Inversion permutes the primitive \(m\)th roots, so \(\Ph\) is reciprocal
and \(u_j=u_{d-j}\). Put \(k=d-1-i\). Reversing the terms after index
\(i\) in \(\Ph(\eta)=0\) gives
\[
   P_i(\eta)=-\eta^d P_k(\eta^{-1}),\qquad
   |P_i(\eta)|=|P_k(\eta)|.
\]
Here \(\eta^{-1}=\overline{\eta}\) and the coefficients of \(P_k\) are
real. Let \(j=\min(i,k)\) and \(B=P_j\). Both \(i\) and \(k\) are
positive, and \(2j\le i+k=d-1\). The product formula
\(\Ph(T)=\prod_{s\mid m}(1-T^s)^{\mu(m/s)}\), with \(\mu\) the
M\"obius function, gives \(u_1=-\mu(m)\ne0\), since only the factor
with \(s=1\) contributes a linear term and \(m\) is squarefree.
Thus \(e=\deg B\) satisfies \(1\le e\le j\), and \(u_e\ne0\)
because \(e\) is the largest index \(t\le j\) for which \(u_t\ne0\).

Suppose that \(C=|P_i(\eta)|^2\) belongs to \(\Q\). By the equality
of moduli above, \(C=|B(\eta)|^2\) for either choice of \(j\). Set
\[
   H(T)=T^e\bigl(B(T)B(T^{-1})-C\bigr).
\]
The exponents in \(B(T)B(T^{-1})\) range from \(-e\) to \(e\), so
\(H\in\Q[T]\) and \(\deg H\le2e<d\). Since \(\eta^{-1}=\overline{\eta}\) and \(B\) has real coefficients,
\(B(\eta^{-1})=\overline{B(\eta)}\). Hence
\[
   H(\eta)=\eta^e\bigl(B(\eta)\overline{B(\eta)}-C\bigr)
           =\eta^e\bigl(|B(\eta)|^2-C\bigr)=0.
\]
The coefficient of \(T^{-e}\) in \(B(T)B(T^{-1})\) is \(u_0u_e\).
Since \(e>0\), it is also the constant coefficient of \(H\), and it is
nonzero. Thus \(H\) is a nonzero polynomial over \(\Q\) of degree
less than \(d\) vanishing at \(\eta\), contrary to the minimality of
\(\Ph\).
\end{proof}

\section{The extremal gaps}
\begin{proof}[Proof of Theorem~\ref{thm:main}]
Set \(F=\Phi_{mp}\), \(D=d(p-1)\), and \(Q=\lfloor p/m\rfloor\).
Proposition~\ref{prop:trunc} gives \(g(F)\le d\). For a gap
\((a,b)\) of length \(d\), take \(i,c,\zeta\) as in \eqref{eq:unit}.
Since \(|P_i(\zeta^p)|^2=c^2\in\Q\), Lemma~\ref{lem:prefix} applied
to \(\eta=\zeta^p\) gives \(i=0\) or \(i=d-1\).

If \(i=0\), then \(1=-c\zeta^b\) and \(0<b\le p\). The only
rational root of unity of odd order is \(1\). Thus \(m\mid b\) and
\(c=-1\). Writing \(b=jm\) gives the first family in \eqref{eq:pairs}.

If \(i=d-1\), then \(P_i(\zeta^p)=-\zeta^{dp}\), so
\(c=\zeta^{dp-b}=1\) and \(m\mid dp-b\).
The inequalities \((d-1)p<b\le D\) give \(d\le dp-b<p\). Write
\(dp-b=jm\). Since \(D=dp-d\), we obtain
\((a,b)=(D-jm,D-jm+d)\), as in the second family.

In both cases \(j\) is a positive integer with \(jm\le p\), hence \(j\le Q\).

\smallskip\noindent\textit{(i) Suppose \(p<m\).}
Then \(Q=0\), so neither case can yield a gap of length \(d\).
By Proposition~\ref{prop:trunc}, all gaps therefore have length less than \(d\).

\smallskip\noindent\textit{(ii) Suppose \(p>m\).}
Then \(Q\ge1\). To realize the candidates in \eqref{eq:pairs}, put \(U(x)=(1-x^m)/\Ph(x)\).
This polynomial has degree \(m-d\), constant coefficient \(1\), and
leading coefficient \(-1\). In \(\Z[[x]]\) we have
\begin{equation}\label{eq:copies}
    \frac{1}{\Ph(x)}=\sum_{h\ge0} x^{hm}U(x).
\end{equation}
By \eqref{eq:product}, \(F(x)\) and \(1/\Ph(x)\) agree in all degrees
less than \(p\). For \(1\le j\le Q\) we have \(jm<p\), since
\(m>1\) and \(p\nmid m\). The \((j-1)\)st translate in
\eqref{eq:copies} ends at exponent \(jm-d\), and the \(j\)th starts
at \(jm\). This gives the first family with endpoint coefficients
\((-1,1)\). The reciprocity of \(F\) gives the second family and
reverses the endpoint coefficients.

The right endpoints in the two families are congruent to \(0\) and
\(dp\) modulo \(m\), respectively. These classes are distinct since
\(\gcd(p,m)=1\) and \(0<d<m\). Each family therefore contributes \(Q\)
distinct gaps, giving \(2Q\) maximum gaps in total.
\end{proof}

\begin{example}\label{ex:odd}
Take \(m=3\cdot5\cdot7=105\) and \(p=211\).
Then \(d=48\),
\(D=10080\), and \(Q=2\). In increasing order the four maximum gaps are
\[
   (57,105),\qquad (162,210),\qquad (9870,9918),\qquad (9975,10023).
\]
\end{example}

\section{Even values of the cofactor}\label{sec:even}
For an even integer \(m\), let \(r\) be the product of its distinct odd
prime divisors, with the empty product equal to \(1\). Write
\[
     m=2hr,\qquad L=m/2=hr,\qquad
     d=\varphi(m)=h\varphi(r).
\]
Throughout this section \(p\) is a prime with \(p\nmid m\), so \(p\) is odd.
Put \(D=d(p-1)\).

\begin{theorem}\label{thm:even}
Let \(m\) be even and let \(p\nmid m\) be a prime. Then
\[
   g(\Phi_{mp})=h\,g(\Phi_{rp})\le d,\qquad
   g(\Phi_{mp})=d\quad\Longleftrightarrow\quad p>r.
\]
If \(r>1\), the gaps of length \(d\) and the coefficients at their
endpoints, in increasing order of exponent, are given by
\[
\begin{array}{c|c}
  \text{gap of length }d & \text{endpoint coefficients} \\ \hline
  (jL-d,jL)     & ((-1)^{j+1},(-1)^j) \\
  (D-jL,D-jL+d) & ((-1)^j,(-1)^{j+1})
\end{array}
\qquad 1\le j\le\lfloor p/r\rfloor .
\]
There are exactly \(2\lfloor p/r\rfloor\) such gaps. If \(r=1\),
all \(p-1\) gaps have length \(d=h\). They are \((kh,(k+1)h)\)
with coefficients \(((-1)^k,(-1)^{k+1})\), for \(0\le k\le p-2\).
\end{theorem}

\begin{proof}
We combine the radical reduction \eqref{eq:radical} with the sign
identity of \cite[\S2]{HLLP}. For odd \(k>1\), negation bijects
primitive \(k\)th roots with primitive
\(2k\)th roots. Since \(\varphi(k)\) is even, this gives
\(\Phi_{2k}(x)=\Phi_k(-x)\). Now \(\operatorname{rad}(mp)=2rp\)
and \(rp>1\), whence
\begin{equation}\label{eq:even-reduction}
    \Phi_{mp}(x)=\Phi_{2rp}(x^h)=\Phi_{rp}(-x^h).
\end{equation}
Thus a nonzero coefficient \(c\) at degree \(a\) in \(\Phi_{rp}\) gives
the coefficient \((-1)^a c\) at degree \(ha\) in \(\Phi_{mp}\).
The ordered supports correspond bijectively, so every gap is scaled by
\(h\) with its multiplicity preserved.

If \(r>1\), apply Theorem~\ref{thm:main} with cofactor \(r\).
For \(p<r\), part~(i) gives strict inequality. For \(p>r\), part~(ii)
gives the maximum gaps of \(\Phi_{rp}\). Writing \(\delta=\varphi(r)\), the two exponents of the first gap
\((jr-\delta,jr)\) both have the parity of \(j\), since \(r\) is odd
and \(\delta\) is even. Multiplication of both coefficients by
\((-1)^j\) gives the first row. Reflection in degree \(D\) reverses
the coefficients and gives the second row. Distinctness and the range
of \(j\) follow from the support bijection. If \(r=1\), the same identity
gives \(\Phi_{mp}(x)=\sum_{k=0}^{p-1}(-1)^k x^{hk}\).
\end{proof}

For squarefree even \(m>2\), one has \(h=1\) and \(r=m/2\). Thus
equality already holds for \(m/2<p<m\), and the count is
\(2\lfloor 2p/m\rfloor\). For \(p<r\), all gaps are multiples of
\(h\), hence \(g(\Phi_{mp})\le d-h\). In that range the actual
maximum is still given by the reduction \eqref{eq:even-reduction}.
For example, if \(r=q\) is prime, apply Theorem~\ref{thm:main}(ii) with
the smaller of \(p,q\) as cofactor. The actual maximum is then
\(h(\min(p,q)-1)\) and occurs
\(2\lfloor\max(p,q)/\min(p,q)\rfloor\) times.

\begin{example}\label{ex:even}
For \(m=30\), \(p=17\), we have \(r=15\), \(d=8\), and \(D=128\).
The extreme pairs are \((7,15)\) and \((113,121)\), with coefficients
\((1,-1)\) and \((-1,1)\), respectively. For \(m=60\) and the same
prime, \(r\) is unchanged and \(h=2\). The extreme pairs become
\((14,30)\) and \((226,242)\) with the same coefficients.
For \(m=30\), \(p=13\), the reduction to \(\Phi_{195}\) instead gives
four maximum gaps of length \(6\) discussed after Theorem~\ref{thm:main}.
\end{example}

\section{Zero runs in the periodic blocks}
The partial-sum argument also identifies the blocks in the coefficient
decomposition of \cite{AAHL} that can contain an extremal run of
zeros. Continue with the notation and hypotheses of
Theorem~\ref{thm:main}, and put \(D=d(p-1)\). Work in
\(R=\Q[x]/(x^m-1)\) and put
\[
     G_i=U(x)P_i(x^p)\ \text{in }R,\qquad 0\le i<d .
\]
For \(G\in R\), let \(G[t]\) denote the coefficient at the residue of
\(t\) modulo \(m\) in the representative of degree less than \(m\).
Multiplication by \(x\) rotates this list of \(m\) coefficients.
For \(0\le t\le w\le D\), the coefficients in a short interval ending
at \(w\) satisfy
\begin{equation}\label{eq:window}
    [x^t]\Phi_{mp}(x)=G_{\lfloor w/p\rfloor}[t]
    \qquad\text{if }w-t\le d-1 .
\end{equation}
Indeed, \eqref{eq:product} and \eqref{eq:copies} give the formula with
\(w=t\). Passing from \(\lfloor t/p\rfloor\) to \(\lfloor w/p\rfloor\)
adds terms \(u_jU[t-jp]\) with \(t<jp\le w\). Each vanishes because
\(1\le jp-t\le d-1\) and \(\deg U=m-d\). Thus \eqref{eq:window}
also holds when the interval crosses several multiples of \(p\).

\begin{proposition}\label{prop:block}
No \(G_i\) contains \(d\) cyclically consecutive zeros. A run of
\(d-1\) cyclically consecutive zeros occurs precisely for \(i=0\) and
\(i=d-1\), once per period in each case. In \(G_0\) it ends at
residue \(-1\), and in \(G_{d-1}\) it ends at residue \(dp-1\).
\end{proposition}
\begin{proof}
Since \(x^m-1=-U\Ph\), cancellation in \(\Q[x]\) gives
\begin{equation}\label{eq:cancel}
    U A=0\text{ in }R\quad\Longleftrightarrow\quad \Ph\mid A
    \qquad(A\in\Q[x]).
\end{equation}
In particular \(G_i\ne0\), because \(P_i(\zeta^p)\ne0\) for
a primitive \(m\)th root \(\zeta\). The polynomial representative of
every rotation of \(G_i\) is divisible by \(U\). Rotate a run of
\(\ell\) zeros to the final \(\ell\) positions. The representative
then has degree at most \(m-1-\ell\). Since \(\deg U=m-d\), this
excludes \(\ell=d\). For \(\ell=d-1\) the representative must equal
\(cU\) for some \(c\in\Q^\times\). Undoing the rotation yields
\(G_i=cx^sU\) in \(R\). By \eqref{eq:cancel},
\(P_i(\zeta^p)=c\zeta^s\). Lemma~\ref{lem:prefix} forces
\(i=0\) or \(i=d-1\).

The last \(d-1\) entries of \(G_0=U\) are zero. If a rotation also
has these final zeros, it equals \(cU\) by the same degree bound.
Equation \eqref{eq:cancel} then gives \(\zeta^s=c\), so \(m\mid s\).
Thus the run in \(U\) is unique. Finally, \eqref{eq:product} and
\eqref{eq:cancel} give \(G_{d-1}=-x^{dp}U\) in \(R\). Its unique
run is therefore shifted by \(dp\) residues.
\end{proof}

\subsection*{Blocks for an even cofactor}
Let \(m\) be even and \(p\nmid m\) prime, and use the notation of
Section~\ref{sec:even}. Define \(U,P_i,G_i\) with
this even \(m\) by the same formulas as above, and use the
representative modulo \(x^m-1\) to read each block.

\begin{proposition}\label{prop:even-block}
Let \(m>2\) be even. No \(G_i\) contains \(d\) cyclically consecutive
zeros. A run of length \(d-1\) occurs precisely for the indices
\[
       0\le i<h\qquad\text{or}\qquad d-h\le i<d.
\]
In each such block there are exactly two runs of length \(d-1\) per
period of \(m\). In the first range they end at residues \(-1,L-1\),
and in the second at \(dp-1,dp+L-1\), all taken modulo \(m\).
For a power of two the two ranges coincide.
\end{proposition}
\begin{proof}
The degree and cancellation argument in Proposition~\ref{prop:block}
applies to every \(m>1\). It excludes \(d\) consecutive zeros, and a
run of length \(d-1\) forces \(G_i=cx^s U\) with \(c\in\Q^\times\).
Hence \(P_i(\zeta^p)=c\zeta^s\) at any primitive \(m\)th root \(\zeta\).

Suppose \(r>1\). As in \eqref{eq:even-reduction} we have
\(\Phi_m(T)=\Phi_r(-T^h)\). Thus the prefix \(P_i(\zeta^p)\) is
the prefix of \(\Phi_r\) through index \(\lfloor i/h\rfloor\),
evaluated at \(-\zeta^{hp}\), a primitive \(r\)th root. Its squared
modulus is \(c^2\). Lemma~\ref{lem:prefix} forces \(\lfloor i/h\rfloor\)
to be \(0\) or \(\varphi(r)-1\), which gives the asserted ranges.
In the first range \(P_i=1\) and \(G_i=U\). In the second
\(P_i=\Phi_m-T^d\), whence \(G_i=-x^{dp}U\) by
\eqref{eq:cancel} and \(\Phi_m(x)\mid\Phi_m(x^p)\).

The last \(d-1\) positions of \(U\) vanish. A rotation with the same
final zeros must equal a rational multiple of \(U\), by the degree
argument. Cancellation implies \(\zeta^s\in\Q\) for its shift \(s\),
so \(s\equiv0,L\pmod m\). The second shift occurs because \(\Phi_m\mid x^L+1\)
gives \(x^L U=-U\) modulo \(x^m-1\). This gives the two listed end
residues for \(U\) and their translates for \(-x^{dp}U\).

Finally, if \(r=1\), then \(m=2h\), \(d=h\), \(\Phi_m=1+x^h\) and
every \(G_i=U=1-x^h\). The two runs end immediately before \(0\)
and \(h\). As \(dp=hp\equiv h\pmod{2h}\), both lists of end
residues agree.
\end{proof}
For \(m=2\) the only block is \(1-x\); both cyclic positions are
nonzero.

\section*{Declaration}
This paper was written by the author in collaboration with GPT 6 Astra.

\end{document}